# A CONVEX STONE-WEIERSTRASS THEOREM & APPLICATIONS

NATHAN S. FELDMAN & PAUL J. MCGUIRE

ABSTRACT. A *convex-polynomial* is a convex combination of the monomials $\{1, x, x^2, \ldots\}$. This paper establishes that the convex-polynomials on $\mathbb{R}$ are dense in $L^p(\mu)$ and weak* dense in $L^\infty(\mu)$, precisely when $\mu([-1, \infty)) = 0$. It is shown that the convex-polynomials are dense in $C(K)$ precisely when $K \cap [-1, \infty) = \emptyset$, where $K$ is a compact subset of the real line. Moreover, the closure of the convex-polynomials on $[-1, b]$ are shown to be the functions that have a convex-power series representation.

A continuous linear operator $T$ on a locally convex space $X$ is convex-cyclic if there is a vector $x \in X$ such that the convex hull of the orbit of $x$ is dense in $X$. The above results characterize which multiplication operators on various real Banach spaces are convex-cyclic. It is shown for certain multiplication operators that every closed invariant convex set is a closed invariant subspace.

## 1. INTRODUCTION

Let $\mathcal{CP}$ denote the convex hull of the set of monomials $\{1, x, x^2, x^3, \ldots\}$ within the vector space of all polynomials. Thus,

$$\mathcal{CP} = \left\{ \sum_{k=0}^n a_k x^k : a_k \geq 0 \text{ for all } 0 \leq k \leq n \text{ and } \sum_{k=0}^n a_k = 1 \right\}.$$

The elements of $\mathcal{CP}$ are called *convex-polynomials*. In this paper we determine when the convex-polynomials are dense in $C(K)$ or in $L^p(\mu)$ for compact sets $K$ and compactly supported measures $\mu$ supported on the real line. In particular we show that $\mathcal{CP}$ is dense in $C([a, b])$ if and only if $b < -1$ and $\mathcal{CP}$ is dense in $L^p(\mu)$ if and only if $\mu\big([-1, \infty)\big) = 0$.

These polynomial approximation problems are closely related to a question about the dynamics of real self-adjoint operators modeled as multiplication operators on $L^2(\mu)$ where $\mu$ is a positive compactly supported regular Borel measure on the real line. This in turn is related to when the closed invariant convex sets for these operators are the same as their closed invariant subspaces. For example, the closed invariant convex sets for the operator $M_x$ of multiplication by $x$ on $L^p(a,b)$ are the same as the closed invariant subspaces of $M_x$ if and only if $b \leq -1$.

The dynamics of a linear operator involves the (predictable or chaotic) behavior of its orbits. If $T$ is a continuous linear operator on a topological vector space $X$ and if $x \in X$, then the *orbit* of $x$ under $T$ is $\mathrm{Orb}(x,T) = \{T^n x\}_{n=0}^{\infty} = \{x, Tx, T^2 x, \ldots\}$. The operator $T$ is called *hypercyclic* if it has an orbit that is dense in $X$ and $T$ is called *cyclic* if it has an orbit whose linear span is dense in $X$.

We are interested in operators that are convex-cyclic (Bermúdez, Bonilla, Feldman [1], León-Saavedra and Romero de la Rosa [7], Rezaei [8]). An operator $T$ on the space $X$ is called *convex-cyclic* if there is a vector $x \in X$ such that the convex hull of its orbit is dense in $X$.

It is known that no normal, subnormal or hyponormal operators are hypercyclic or supercyclic (Bourdon [2]). However, there are normal operators that are convex-cyclic, in fact, a diagonal matrix on $\mathbb{C}^n$ is convex-cyclic if and only if its diagonal entries $\{\lambda_k\}_{k=1}^n$ are distinct, belong to the set $\mathbb{C} \setminus (\mathbb{R} \cup \overline{\mathbb{D}})$ and satisfy $\lambda_j \neq \overline{\lambda_k}$ for all $1 \leq j, k \leq n$ (see [6]). Similarly, a real diagonal matrix acting on $\mathbb{R}^n$ is convex-cyclic if and only if its diagonal entries are distinct and all belong to the interval $(-\infty, -1)$. The same results are also true for infinite diagonal matrices acting on $\ell_{\mathbb{C}}^2$ or $\ell_{\mathbb{R}}^2$ (see [1]). Thus, there are normal convex-cyclic operators in finite and infinite dimensions. We consider when multiplication operators on certain Banach spaces are convex-cyclic. For example, we show that the operator $M_x$ on $L^p(a,b)$ is convex-cyclic if and only if $b \leq -1$.

In the next section we use the Hahn-Banach Separation Theorem to introduce the basic tool that is used to determine when a convex set is dense in a given space.

## 2. Preliminaries

Recall that $\mathcal{CP}$ denotes the convex hull of the set of monomials $\{1, x, x^2, x^3, \ldots\}$ and the elements of $\mathcal{CP}$ are called *convex-polynomials*. If $f$ has a power series that converges near zero, say $f(x) = \sum_{k=0}^{\infty} a_k x^k$, then the coefficients of the series are related to the derivatives of $f$ by $a_k = f^{(k)}(0)/k!$ for $k \geq 0$. This can be used to easily characterize the convex-polynomials and determine some basic properties of them, as stated in the following proposition.

**Proposition 2.1.** *(Properties of Convex-Polynomials) The following hold:*

(1) *A polynomial $p(x)$ is a convex-polynomial if and only if $p^{(k)}(0) \geq 0$ for all $k \geq 0$ and $p(1) = 1$.*
(2) *$|p(x)| \leq p(|x|)$ for all $x \in \mathbb{R}$ and all $p \in \mathcal{CP}$.*
(3) *$p(\mathbb{R}) \subseteq \mathbb{R}$, $p(\mathbb{R}^+) \subseteq \mathbb{R}^+$ and $p([-1,1]) \subseteq [-1,1]$ for all $p \in \mathcal{CP}$.*
(4) *The set $\mathcal{CP}$ of convex-polynomials is closed under multiplication and composition.*

**Theorem 2.2.** *(Hahn-Banach Criterion for a Convex Set to be Dense) If $C$ is a convex set in a real locally convex space $X$, then $C$ is dense in $X$ if and only if for every non-zero continuous linear functional $f$ on $X$ we have that $\sup_{x \in C} f(x) = \infty$. Furthermore if $C$ is the convex-hull of a set $S$, then $\sup_{x \in C} f(x) = \sup_{x \in S} f(x)$.*

The previous result is a simple restatement of the geometric form of the Hahn-Banach Theorem which says that whenever a point does not belong to a closed convex set, then the point and the convex set can be strictly separated by a real hyperplane. See [5, Theorem 3.13, p. 111].

**Theorem 2.3.** *(Hahn-Banach Characterization of Convex-Cyclicity) Let $X$ be a real locally convex space, $T : X \to X$ a continuous linear operator, and $x \in X$. Then the following are equivalent:*
  (1) *The convex hull of the orbit of $x$ under $T$ is dense in $X$.*
  (2) *For every non-zero continuous linear functional $f$ on $X$ we have $\sup f(Orb(x,T)) = \infty$.*

*Proof.* Apply the Hahn-Banach Separation Theorem. □

If $S \subseteq \mathbb{R}$, then let $M_c(S)$ be the space of all real compactly supported regular Borel measures supported on $S$ with finite total variation. That is, $\mu \in M_c(S)$ if and only if there is a compact set $K$ such that $K \subseteq S$ and $|\mu|(\mathbb{R} \setminus K) = 0$ and $|\mu|(K)| < \infty$. Also let $M_c^+(S)$ be the subspace of $M_c(S)$ consisting of positive measures.

The following result is a corollary of Theorem 2.3 where $T$ is the linear operator of multiplication by the independent variable $x$ on an $L^p$ space.

**Corollary 2.4.** *(Hahn-Banach Characterization) The following hold:*
  (1) *If $\mu \in M_c^+(\mathbb{R})$ and $1 \leq p < \infty$, then the convex-polynomials are dense in $L^p(\mu)$ if and only if for every nonzero $f \in L^q(\mu)$, where $\frac{1}{p} + \frac{1}{q} = 1$, we have $\sup_{n \geq 1} \int x^n f(x) d\mu = \infty$.*
  (2) *If $\mu \in M_c^+(\mathbb{R})$, then the convex-polynomials are weak*-dense in $L^\infty(\mu)$ if and only if for every nonzero $f \in L^1(\mu)$ we have $\sup_{n \geq 1} \int x^n f(x) d\mu = \infty$.*
  (3) *If $K$ is a compact set in the plane, then the convex-polynomials are dense in $C(K)$ if and only if for every nonzero measure $\nu \in M(K)$ we have $\sup_{n \geq 1} \int x^n d\nu = \infty$.*

We are interested in when the convex-polynomials are dense in $L^p(\mu)$. Next we state some natural necessary conditions. We shall see that the necessary condition in case $(a)$ below is also a sufficient condition. However, that is not true in case $(b)$.

**Proposition 2.5.** *(Necessary Conditions) If $\mu \in M_c^+(\mathbb{R})$ and the convex-polynomials are dense in $L^p(\mu)$, then $\mu([-1, \infty)) = 0$.*

The previous result follows easily from Proposition 2.1.

3. The Convex-Polynomials are weak*-dense in $L^\infty(a, -1)$

**Theorem 3.1.** *If $a < b \leq -1$, then the set of convex-polynomials is weak*-dense in $L^\infty(a,b)$.*

Since the weak* dual of $L^\infty(a,b)$ is $L^1(a,b)$, Corollary 2.4 shows that Theorem 3.1 will follow immediately from Theorem 3.2.

**Theorem 3.2.** *(Unbounded Moments) If $a < b \leq -1$ and $f \in L^1(a,b)$ and $f$ is non-zero on a set of positive Lebesgue measure, then*

$$\sup_{n \geq 0} \int_a^b x^n f(x)dx = \infty \ \text{and} \ \inf_{n \geq 0} \int_a^b x^n f(x)dx = -\infty.$$

Before we prove Theorem 3.2 we must establish several preliminary results. If $x_0 \in [a,b]$, then we will say that $x_0$ is a *peak point* for the set $\mathcal{CP}$ of convex-polynomials on $[a,b]$ if there exists a $p \in \mathcal{CP}$ such that $|p(x_0)| > |p(x)|$ for all $x \in [a,b] \setminus \{x_0\}$.

**Theorem 3.3.** *(Peak Points of $\mathcal{CP}$) If $a < -1$, then every point $x_0 \in [a,-1]$ is a peak point for the set of convex-polynomials on $[a,-1]$.*

*Proof.* If $x_0 = -1$, then $p(x) = \frac{1}{1-a}(x-a)$ is a convex polynomial which peaks on $[a,-1]$ at $x_0$. The case when $a \leq x_0 < -1$ is covered by the following proposition. $\square$

**Proposition 3.4.** *(Peaking Convex-Polynomials) If $a < -1$ and $x_0 \in [a,-1)$, then there exists a positive even integer $n$ and an $\alpha \in \left(0, \frac{1}{1-\frac{n+1}{n}x_0}\right)$ such that the polynomial*

$$p_{\alpha,n}(x) = \left(1 - \alpha + \frac{n+1}{n}\alpha x_0\right) - \frac{n+1}{n}\alpha x_0 x^n + \alpha x^{n+1}$$

*has the following properties:*

(1) $p_{\alpha,n}$ is a convex-polynomial;
(2) $p_{\alpha,n}$ has its only critical points at $x_0$ and $0$;
(3) $p_{\alpha,n}$ is strictly increasing on $(-\infty, x_0]$ and strictly decreasing on $[x_0, 0]$;
(4) The absolute maximum of $p_{\alpha,n}$ on $(-\infty, 0]$ is at $x_0$, moreover,

$$0 < p_{\alpha,n}(x) < p_{\alpha,n}(x_0) \ \text{for all} \ x \in [a, 0] \setminus \{x_0\};$$

(5) $p_{\alpha,n}$ is concave down on $(-\infty, \frac{n-1}{n}x_0]$ and concave up on $[\frac{n-1}{n}x_0, 0]$;
(6) 
$$p_{\alpha,n}(x_0) = 1 - \frac{\alpha}{n}\left(x_0^{n+1} - (n+1)x_0 + n\right)$$
$$= 1 - \frac{\alpha}{n}(x_0 - 1)^2 \left(x_0^{n-1} + 2x_0^{n-2} + 3x_0^{n-3} + \cdots + (n-1)x_0 + n\right);$$

(7) For $x_0^{n+1} - (n+1)x_0 + n < 0$, $p_{\alpha,n}(x_0) > 1$.

*Proof.* (1) Clearly the three coefficients of $p_{\alpha,n}$ sum to one, so it suffices to show that they are all positive. Since $\alpha > 0$ and $x_0 < 0$, it follows that the coefficients of $x^n$ and $x^{n+1}$ are both positive. Also, since $\alpha > 0$ we have $\alpha \in \left(0, \frac{1}{1-\frac{n+1}{n}x_0}\right)$ if and only if $\left(1 - \alpha + \frac{n+1}{n}\alpha x_0\right) > 0$. Thus $p_{\alpha,n}$ is a convex-polynomial.

Using the fact that $n$ is even, it is a straightforward calculus exercise to verify (2), (3), (5), and (6). While (7) follows from (6), it should be noted that as $x_0$ nears $-1$, the value of $n$ must increase in order that $(*)$ $x_0^{n+1} - (n+1)x_0 + n < 0$ holds. But once $x_0 < -1$ has been given, we will choose $n$ (first, before $\alpha$) such that $(*)$ holds. Then since $a^{n+1}$ is the dominant term in

$$p_{\alpha,n}(a) = \left(1 - \alpha + \frac{n+1}{n}\alpha x_0\right) - \alpha a^{n+1}\left(\frac{n+1}{n}\frac{x_0}{a} - 1\right),$$

choosing $\alpha$ sufficiently small ($n$ has already been chosen) will result in $p_{\alpha,n}(a) > 0$. Since $p_{\alpha,n}(a)$ is the absolute minimum of $p_{\alpha,n}(x)$ on $[a, x_0]$ and $p_{\alpha,n}(0)$ (which is positive) is the absolute minimum of $p_{\alpha,n}(x)$ on $[x_0, 0]$, (4) follows. □

**Proposition 3.5.** *If $a < b \leq -1$ and $f \in L^1(a,b)$ and $f(x) \geq \varepsilon > 0$ for a.e. $x \in (a,b)$, then the following hold:*

(1) $\int_a^b x^{2n} f(x) dx \to \infty$ *as* $n \to \infty$.
(2) $\int_a^b x^{2n+1} f(x) dx \to -\infty$ *as* $n \to \infty$.

*Proof.* (1) Since $f(x) \geq \varepsilon$ on $(a,b)$ and since $x^{2n} \geq 0$, then $x^{2n} f(x) \geq \varepsilon x^{2n}$ on $(a,b)$, thus we have that

$$\int_a^b x^{2n} f(x) dx \geq \varepsilon \int_a^b x^{2n} dx = \frac{\varepsilon(b^{2n+1} - a^{2n+1})}{2n+1} = \frac{\varepsilon(|a|^{2n+1} - |b|^{2n+1})}{2n+1} \to \infty$$

as $n \to \infty$.

(2) In this case we have $f(x) \geq \varepsilon$ on $(a,b)$ and since $b \leq -1$, $x^{2n+1} \leq 0$ on $(a,b)$, thus $x^{2n+1} f(x) \leq \varepsilon x^{2n+1}$ on $(a,b)$, so we have that

$$\int_a^b x^{2n+1} f(x) dx \leq \varepsilon \int_a^b x^{2n+1} dx = \frac{\varepsilon}{2n+2}(b^{2n+2} - a^{2n+2}) \to -\infty \text{ as } n \to \infty,$$

since $a < b < 0$ and thus $|a| > |b|$. □

**Lemma 3.6.** *If $\mu \in M_c^+(\mathbb{R})$ and $f \in L^1(\mu)$, then the following hold:*

(1) $\sup_{n \geq 1} \int x^n f(x) d\mu = \infty \Leftrightarrow \sup_{p \in \mathcal{CP}} \int p(x) f(x) d\mu = \infty$.

(2) *If $m > 1$ and $c > 0$, then there is a sequence of positive integers $\{n_k\}$ such that*

$$\int x^{n_k} f(x) d\mu \geq cm^k \text{ for all } k \geq 1$$

*if and only if there is a sequence of convex-polynomials $\{p_k\}$ such that*

$$\int p_k(x) f(x) d\mu \geq cm^k \text{ for all } k \geq 1.$$

*In this case, we may choose the $n_k$'s to satisfy $1 \leq n_k \leq \deg(p_k)$ for all $k$.*

*Proof.* (1) The forward implication ($\Rightarrow$) is clearly true. We will prove the converse ($\Leftarrow$) by proving its contrapositive. So, suppose that there exists an $M > 0$ such that $\sup_{n \geq 1} \int x^n f(x) d\mu \leq M$ for all $n \geq 0$ and we will show that $\int p(x) f(x) d\mu \leq M$ for all convex-polynomials $p$. By assumption, the moments $\int x^n f(x) d\mu$ all lie in the convex set $\{x : x \leq M\}$ and since $\int p(x) f(x) d\mu$ is a convex-combination of a finite number of these moments, then $\int p(x) f(x) d\mu$ also belongs to the convex set $\{x : x \leq M\}$.

(2) The forward implication ($\Rightarrow$) is clearly true. We will prove the reverse direction ($\Leftarrow$). So suppose that there exists a sequence of convex-polynomials $\{p_k\}$ such that for every $k \geq 1$

$$(*) \qquad \int p_k(x) f(x) d\mu \geq cm^k.$$

then we will construct a sequence $\{n_k\}$ with the required property. Fix a $k \geq 1$, then we have that $\int p_k(x) f(x) d\mu(x) \geq cm^k$. Suppose that $p_k(x) = \sum_{j=1}^{M} a_j x^j$ where $a_j \geq 0$ and $\sum_{j=1}^{M} a_j = 1$. Since $\int p_k(x) f(x) d\mu$ is a convex combination of the values $\{\int x^j f(x) d\mu : 1 \leq j \leq M\}$, if all of these values are in the convex set $H := \{x : x < cm^k\}$, then a convex combination of them would also belong to $H$, contradicting equation $(*)$ above. Thus there is some $n_k \in \{0, 1, 2, \ldots, deg(p_k)\}$ such that $\int x^{n_k} f(x) d\mu$ does not belong to $H$, thus $\int x^{n_k} f(x) d\mu \geq cm^k$. This produces a sequence $\{n_k\}$ with the required property. $\square$

**Proposition 3.7.** *Suppose that $a < b \leq -1$ and that $f \in L^1(a,b)$. If $f$ is positive on an open interval $I_0 \subseteq (a,b)$, then there is an $m > 1$, a number $c > 0$, and a sequence of positive integers $n_k \to \infty$ such that*

$$\int_a^b x^{n_k} f(x) dx \geq cm^k \to \infty \text{ as } k \to \infty.$$

*Proof.* By Lemma 3.6 it suffices to find an $m > 1$, a $c > 0$, and a sequence of convex-polynomials $\{p_k\}_{k=1}^{\infty}$ such that

$$\int_a^b p_k(x) f(x) dx \geq cm^k \text{ for all } k \geq 1.$$

Let $I_0 := (e, d) \subseteq (a, b)$ be an open interval on which $f$ is positive and let $x_0 \in I_0$. Since $b \leq -1$, then $x_0 < -1$, thus by Proposition 3.4, there exists a convex-polynomial $p$ that peaks on $[a, b]$ at $x_0$ and such that $p(x_0) > 1$. By the continuity of $p$ at $x_0$ and the fact that $x_0$ is the unique point in $[a, b]$ at which the absolute maximum of $p$ is attained, we can choose an $m$ such that $1 < m < p(x_0)$ and such that if $p(x) \geq m$ for some $x \in [a, b]$, then $x \in I_0$. Now let $I_1 = [a, b] \cap p^{-1}([m, p(x_0)])$. Then $I_1$ is a closed interval and $x_0 \in I_1 \subseteq I_0 \subseteq (a, b)$ and $p(x) \geq m$ for all $x \in I_1$. Thus we have the following:

$$\int_a^b p(x)^k f(x) dx = m^k \left[ \int_{I_1} \left( \frac{p(x)}{m} \right)^k f(x) dx + \int_{[a,b] - I_1} \left( \frac{p(x)}{m} \right)^k f(x) dx \right].$$

If $x \in [a,b] - I_1$, then $0 \leq p(x) < m$, so $0 \leq \frac{p(x)}{m} < 1$, thus $\left(\frac{p(x)}{m}\right)^k \to 0$ as $k \to \infty$, so

$$\delta_k := \int_{[a,b]-I_1} \left(\frac{p(x)}{m}\right)^k f(x)dx \to 0$$

as $n \to \infty$ by the dominated convergence theorem and the fact that $f \in L^1(a,b)$.

If $x \in I_1$, then $p(x) \geq m$, so $(\frac{p(x)}{m})^k \geq 1$ for all $k \geq 1$. Now since $f(x) > 0$ for $x \in I_1$, we can multiply the previous inequality by $f$ and preserve the inequality, giving $(\frac{p(x)}{m})^k f(x) \geq f(x)$ for all $x \in I_1$. This then gives,

$$\int_{I_1} \left(\frac{p(x)}{m}\right)^k f(x)dx \geq \int_{I_1} f(x)dx =: \varepsilon_0 > 0.$$

It follows that for large $k$ we have,

$$\int_a^b p(x)^k f(x)dx \geq m^k[\varepsilon_0 + \delta_k] \geq m^k[\varepsilon_0 - \frac{1}{2}\varepsilon_0] = \frac{1}{2}\varepsilon_0 \cdot m^k \to \infty$$

since $k \to \infty$, $\varepsilon_0 > 0$, and $m > 1$. The result now follows by Lemma 3.6. □

We are now ready for the proof of Theorem 3.2.

*Proof.* (Proof of Theorem 3.2) Let $a < -1$. In order to show that the convex-polynomials are weak* dense in $L^\infty(a,-1)$, it suffices, by Corollary 2.4, to show that for every non-zero function $f$ in $L^1(a,-1)$ that

(*) $$\sup_{n \geq 1} \int_a^{-1} x^n f(x)dx = \infty.$$

So let $f$ be a non-zero function in $L^1(a,-1)$. For $x \in [a,-1]$, let $F(x) = \int_a^x f(t)dt$. Since $f \in L^1(a,-1)$, the function $F$ is an absolutely continuous non-zero function on $[a,b]$. As such we may use integration by parts on the above integral which gives

(**) $$\int_a^{-1} x^n f(x)dx = x^n F(x)|_a^{-1} - \int_a^{-1} nx^{n-1}F(x)dx =$$

$$= [(-1)^n F(-1) - 0] - \int_a^{-1} nx^{n-1}F(x)dx.$$

Now consider two cases, in each case we will show that $\sup_{n \geq 1} \int_a^{-1} x^n f(x)dx = \infty$.

*Case 1:*  There is an $x_0 \in [a,-1]$ such that $F(x_0) < 0$.

Referring to equation (**) above we will show that

$$-\int_a^{-1} nx^{n-1}F(x)dx = \int_a^{-1} nx^{n-1}(-1)F(x)dx$$

has a subsequence that converges to infinity. Since $(-1)F(x_0) > 0$ and $F$ is continuous on $[a,-1]$, then there is an open interval contained in $[a,-1]$ on which $(-1)F$ is positive. Thus we may apply Proposition 3.7 which says that there is a sequence $n_k \to \infty$ such that $\int_a^{-1} x^{n_k}(-1)F(x)dx \to \infty$. Letting $n = n_k + 1$ in equation $(**)$ above gives the following:

$$\int_a^{-1} x^{n_k+1} f(x)dx = (-1)^{n_k+1} F(-1) + (n_k + 1) \int_a^{-1} x^{n_k}(-1) F(x)dx \to \infty.$$

Thus $(*)$ is satisfied in this case.

*Case 2:* $F(x) \geq 0$ for all $x \in [a, -1]$.

In this case we may apply Proposition 3.5 to see that $\int_a^{-1} x^{2n+1} F(x) dx \to -\infty$ as $n \to \infty$ and thus, referring to $(**)$ we have

$$\int_a^{-1} x^{2n+2} f(x)dx = x^{2n+2} F(x)|_a^{-1} - \int_a^{-1} (2n+2) x^{2n+1} F(x) dx =$$

$$= (1)F(-1) - 0 - (2n+2) \int_a^{-1} x^{2n+1} F(x)dx \to +\infty.$$

Thus, $(*)$ is satisfied in this case. Since $(*)$ holds in both cases, the result follows. □

**Corollary 3.8.** *If $1 \leq p < \infty$ and $a < -1$, then the convex-polynomials are dense in $L^p(a, -1)$.*

*Proof.* From Theorem 3.1 we get that the convex-polynomials are weak*-dense in the space $L^\infty(a, -1)$, and since $L^\infty(a, -1)$ is norm dense in $L^p(a, -1)$ when $1 \leq p < \infty$, the result follows. □

4. AN INTERLUDE ON RIEMANN-STIELTJES INTEGRALS

This is a brief review of Riemann-Stieltjes integrals and their relation to Lebesgue integrals. This is all classic material that can be found in many souces. A modern reference is [4, p. 75-84]. Suppose that $g$ and $F$ are bounded functions on an interval $[a, b]$ and $P = \{x_0, x_1, \ldots, x_n\}$ is a *partition* of $[a, b]$, meaning $a = x_0 < x_1 < \ldots < x_n = b$, and let $c_k \in [x_{k-1}, x_k]$ for each $1 \leq k \leq n$. Then the *Riemann-Stieltjes sum* $S(g, F, P, \{c_k\})$ is defined to be

$$S(g, F, P, \{c_k\}) = \sum_{k=1}^n g(c_k)[F(x_k) - F(x_{k-1})].$$

By a *tagged partition* of $[a, b]$ we mean a partition $P = \{x_0, \ldots, x_n\}$ of $[a, b]$ together with a choice of points (called tags) $c_k \in [x_{k-1}, x_k]$ for $1 \leq k \leq n$. Also, the norm of the partition $P$ is defined to be $\|P\| := \max\{(x_k - x_{k-1}) : 1 \leq k \leq n\}$.

A function $g$ defined on $[a, b]$ is *Riemann-Stieltjes integrable* with respect to $F$ if there exists a number $A$ such that for every $\varepsilon > 0$, there is a $\delta > 0$ such that for every tagged partition $P$ with norm less than $\delta$ we have that

$$|S(g, F, P, \{c_k\}) - A| < \varepsilon.$$

The number $A$, when it exists, is called the Riemann-Stieltjes integral of $g$ with respect to $F$ and is denoted by $\int_a^b g(x) dF(x)$.

**Proposition 4.1.** *(Existence of the Riemann-Stieltjes Integral) If $g$ is continuous on $[a,b]$ and $F$ is a function of bounded variation (a difference of two monotone increasing functions) on $[a,b]$, then $\int_a^b g(x)dF(x)$ exists.*

Continuing from the above proposition, it can also be shown that if we define $G(x) = \int_a^x g(t)dF(t)$, then $G$ is continuous wherever $F$ is continuous and $G$ is differentiable wherever $F$ is differentiable and at such a point $x$ we have $G'(x) = g(x)F'(x)$. Recall that functions of bounded variation are differentiable almost everywhere with respect to Lebesgue measure. The proofs of the above remarks, the previous proposition and the following proposition can be found in [4, p. 79-82].

**Proposition 4.2.** *(Integration by Parts) If $g$ and $F$ are bounded functions on $[a,b]$ with no common discontinuities in $[a,b]$ and if $\int_a^b g(x)dF(x)$ exists, then $\int_a^b F(x)dg(x)$ exists and satisfies $\int_a^b g(x)dF(x) = g(b)F(b) - g(a)F(a) - \int_a^b F(x)dg(x)$.*

Now if $\mu$ is a finite regular Borel measure on $\mathbb{R}$, then its distribution function $F : \mathbb{R} \to \mathbb{R}$ is defined by $F(x) = \mu((-\infty, x])$. It follows that $F$ is a bounded function of bounded variation and hence a difference of two monotone increasing functions. In particular, $F$ is continuous on $\mathbb{R}$ except at a countable set of points.

**Proposition 4.3.** *(Lebesgue & Riemann-Stieltjes Integrals) If $\mu$ is a finite regular Borel measure on $\mathbb{R}$, $F(x) = \mu((-\infty, x])$ is its distribution function and $g$ is a continuous function on a closed interval $[a,b]$, then both the Lebesgue integral $\int_{[a,b]} g(x)d\mu$ and the Riemann-Stieltjes integral $\int_a^b g(x)dF(x)$ exist and they are equal: $\int_{[a,b]} g(x)d\mu = \int_a^b g(x)dF(x)$.*

*Proof.* Since $g$ is continuous on $[a,b]$, it is known that both integrals exist. To see that they are equal we will show that for every $\varepsilon > 0$, there is a Riemann-Stieltjes sum that is within $\varepsilon$ of each integral. This implies that they must be equal. So let $\varepsilon > 0$, since $\int_a^b g(x)dF(x)$ exists, there is a $\delta > 0$ such that for any tagged partition $P$ with $\|P\| < \delta$ we have that the Riemann-Stieltjes sum $S(g, F, P, \{c_k\})$ is within $\varepsilon$ of $\int_a^b g(x)dF(x)$. Also, since $g$ is continuous on $[a,b]$, then it is uniformly continuous on $[a,b]$, so it may be uniformly approximated by a step function of the form $\varphi(x) = \sum_{k=1}^n g(c_k)\chi_{(x_{k-1}, x_k]}$ where the partition $P = \{x_0, x_1, \ldots, x_n\}$ has norm at most $\delta$ and so that $\int_{[a,b]} \varphi(x)d\mu$ is within $\varepsilon$ of $\int_{[a,b]} g(x)d\mu$. Then simply notice that $\int_{[a,b]} \varphi(x)d\mu$ is a Riemann-Stieltjes sum, in fact it is equal to $\sum_{k=1}^n g(c_k)[F(x_k) - F(x_{k-1})]$. Since $\|P\| < \delta$, then this sum is also within $\varepsilon$ of $\int_a^b g(x)dF(x)$. So, this sum is within $\varepsilon$ of each integral. Since $\varepsilon > 0$ is arbitrary, the two integrals must be equal. □

## 5. When are the Convex-Polynomials dense in $L^p(\mu)$?

**Theorem 5.1.** *(Moments of Measures) If $\mu$ is a (positive or signed) finite real measure with compact support in $(-\infty, -1]$ that is not supported on the singleton*

*set $\{-1\}$, then*

$$\sup_{n \geq 0} \int x^n d\mu = \infty \text{ and } \inf_{n \geq 0} \int x^n d\mu = -\infty.$$

*Proof.* Let $F : \mathbb{R} \to \mathbb{R}$ be the distribution function for $\mu$ given by $F(x) = \mu((-\infty, x])$. Since $\mu$ has compact support we may choose an $a < -1$ such that $\mu$ is the zero measure on $(-\infty, a]$. Thus $F(a) = 0$. Notice that since $\mu$ is not supported on $\{-1\}$, then $\mu$ is not the zero measure on $(a, -1)$, thus $F$ is not almost everywhere equal to zero on $(a, -1)$. Since $x^n$ is continuous we may use integration by parts (Proposition 4.2) and a Riemann-Stieltjes integral (Proposition 4.3) to give

$$\int x^n d\mu = \int_{[a,-1]} x^n d\mu =$$

$$\int_{\{-1,a\}} x^n d\mu + \int_{(a,-1)} x^n d\mu = (-1)^n \mu(\{-1\}) + a^n \mu(\{a\}) + \int_{(a,-1)} x^n d\mu =$$

$$(-1)^n \mu(\{-1\}) + a^n \mu(\{a\}) + \int_a^{-1} x^n dF(x) =$$

$$(-1)^n \mu(\{-1\}) + a^n \cdot 0 + \left[ x^n F(x) \big|_a^{-1} - \int_a^{-1} n x^{n-1} F(x) dx \right] =$$

$$(-1)^n \mu(\{-1\}) + 0 + \left[ (-1)^n F(-1) - 0 - \int_a^{-1} n x^{n-1} F(x) dx \right].$$

Since $F \in L^\infty(a, -1)$ and $F$ is non-zero on a set of positive Lebesgue measure in $(a, -1)$, then by Theorem 3.2 we have that $\inf_{n \geq 1} \int_a^{-1} n x^{n-1} F(x) dx = -\infty$, so $\sup_{n \geq 1} (-1) \int_a^{-1} n x^{n-1} F(x) dx = \infty$. It follows that $\sup_{n \geq 0} \int x^n d\mu = \infty$. Since the previous statement holds for all measures with compact support in $(-\infty, -1]$ that are not supported on $\{-1\}$, then we may apply it to the measure $(-1)\mu$ and also conclude that $\inf_{n \geq 0} \int x^n d\mu = -\infty$. □

In the result below recall that we are considering real-valued $L^p$ functions.

**Corollary 5.2.** *For a finite positive regular Borel measure $\mu$ with compact support on the real line, the following are equivalent:*

(1) *The convex-polynomials are weak\*-dense in $L^\infty(\mu)$;*
(2) *The convex polynomials are dense in $L^p(\mu)$ for all $1 \leq p < \infty$;*
(3) *For any $f \in L^p(\mu)$ with $|f| > 0$ $\mu$ a.e. the set $\{p(x) \cdot f(x) : p \in \mathcal{CP}\}$ is dense in $L^p(\mu)$ when $1 \leq p < \infty$ and weak\*-dense in $L^\infty(\mu)$ when $p = \infty$;*
(4) $\mu([-1, \infty)) = 0$.

*Proof.* By Proposition 2.1 convex-polynomials are all positive on $(0, \infty)$ and they are all bounded by one in absolute value on the interval $[-1, 1]$, thus it follows that conditions $(1), (2)$, and $(3)$ each imply condition $(4)$. Conversely, if $(4)$ holds, then using the Hahn-Banach characterization of when a convex set is dense (Theorem 2.2) together with Theorem 5.1, it follows that conditions $(1), (2)$, and $(3)$ each hold. □

## 6. Stone-Weierstrass Theorems for Convex-Polnomials

6.1. **Topologies on Spaces of Continuous Functions.** If $X$ is a locally compact Hausdorff topological space, then let $C(X)$ denote the set of all continuous functions on $X$, let $C_b(X)$ denote the set of all bounded continuous functions on $X$, and let $C_0(X)$ denote the set of all continuous functions on $X$ that "vanish at infinity"; that is $f \in C_0(X)$ if $f$ is continuous on $X$ and for every $\varepsilon > 0$, there is a compact set $K$ in $X$ such that $|f| < \varepsilon$ on $X \setminus K$. If $X$ is an open interval $(a,b)$ on the real line, then $C_0(a,b)$ is the space of all continuous functions on $(a,b)$ which "vanish at the endpoints"; meaning $\lim_{x \to a^+} f(x) = \lim_{x \to b^-} f(x) = 0$.

If $C(X)$ is endowed with the topology of uniform convergence on compact sets, then its dual space is $M_c(X)$, the space of all regular Borel measures on $X$ with compact support in $X$, see [5, p. 114]. Following Buck [3], define the *strict topology* on the space $C_b(X)$ to be the topology given by the family of semi-norms

$$\|f\|_\varphi = \|\varphi f\|_\infty = \sup_{x \in X} |\varphi(x) f(x)|$$

where $\varphi \in C_0(X)$. In [3, Theorem 2] it is shown that the dual space of $C_b(X)$ endowed with the strict topology is the space $M(X)$ of all finite regular Borel measures on $X$. In what follows we will take $X$ to be an open interval $(a,b)$ on the real line.

**Theorem 6.1.** *If $a < b \leq -1$, then the set of convex-polynomials is dense in $C_b(a,b)$ with the strict topology.*

*Proof.* This follows from the Hahn-Banach Criterion (Theorem 2.2) for the density of a convex set. Since the dual space of $C_b(a,b)$ is the set $M(a,b)$ (see [3, Theorem 2]), and since $b \leq -1$, all such non-zero measures have unbounded moments by Theorem 5.1. □

**Theorem 6.2.** *If $a < b \leq -1$, then the set of convex-polynomials is dense in $C(a,b)$ with the topology of uniform convergence on compact subsets of $(a,b)$.*

*Proof.* This follows from the Hahn-Banach Criterion (Theorem 2.2) for the density of a convex set. Since the dual space of $C(a,b)$ is the set $M_c(a,b)$, (see [5, p. 114]), and since $b \leq -1$, all such non-zero measures have unbounded moments by Theorem 5.1. □

**Theorem 6.3.** *(A Convex Stone-Weierstrass Theorem) If $a < b$, then the convex-polynomials are uniformly dense in $C([a,b])$ if and only if $b < -1$.*

*Proof.* If $b < -1$, then the density of the convex-polynomials in $C[a,b]$ follows from the Hahn-Banach Criterion (Theorem 2.2) for the density of a convex set. Since the dual space of $C[a,b]$ is the set $M([a,b])$ of all regular Borel measures with finite total variation (see [5, p. 75]), and since $b \leq -1$, all such non-zero measures have unbounded moments by Theorem 5.1.

For the converse, simply notice that all convex-polynomials are positive on $(0, \infty)$ and they are all bounded by one in absolute value on the interval $[-1, 1]$ (by Proposition 2.1), thus they cannot approximate the function $f(x) = -5$ on $[a, b]$ when $b \geq -1$. □

**Corollary 6.4.** *Suppose that $a < b < -1$ and $k$ is a positive integer. Then the following hold:*

*(a) the convex-hull of the set $\{x^{kn} : n \geq 0\}$ is dense in $C(a, -1)$ with the topology of uniform convergence on compact sets if and only $k$ is odd.*

*(b) the convex-hull of the set $\{x^{kn} : n \geq 0\}$ is dense in $C[a, b]$ with the topology of uniform convergence on $[a, b]$ if and only $k$ is odd.*

*Proof.* (a) Let $a < -1$ and suppose that $k$ is an odd positive integer and we will show that the set $\{p(x^k) : p \in \mathcal{CP}\}$ is dense in $C(a, -1)$. Let $g : (a, -1) \to (a^k, -1)$ be given by $g(x) = x^k$. Notice that since $k$ is odd, that $a^k < (-1)^k = -1$. Also notice that the operator $C_{x^k}$ of composition with $x^k$, maps $C(a^k, -1) \to C(a, -1)$ and composition with $x^{1/k}$, $C_{x^{1/k}} = C_{x^k}^{-1} : C(a, -1) \to C(a^k, -1)$. Let $f \in C(a, -1)$. Then $f(x^{1/k}) \in C(a^k, -1)$. By Theorem 6.2, the convex-polynomials are dense in $C(a^k, -1)$, so there exists a sequence $\{p_j\}$ of convex-polynomials such that $p_j(x) \to f(x^{1/k})$ uniformly on compact subsets of $(a^k, -1)$. It then follows that $p_j(x^k) \to f(x)$ uniformly on compact subsets of $(a, -1)$. Since $f \in C(a, -1)$ was arbitrary, we see that the set $\{p(x^k) : p \in \mathcal{CP}\}$ is dense in $C(a, -1)$.

For the converse, simply notice that if $k$ is even, then $x^{kn} \geq 0$ for all $x$ and all $n$, thus their convex-hull cannot be dense in $C(a, -1)$.

(b) This follows, as in the proof of Theorem 6.3, simply because if $a_1 < a$, then $[a, b]$ is a compact subset of $(a_1, -1)$ and continuous functions on $[a, b]$ can be extended to be continuous on $(a_1, -1)$, thus the result follows from part $(a)$. □

**Corollary 6.5.** *If $a < -1$ and $k$ is an odd positive integer, then for every nonzero $\mu \in M((a, -1))$ we have*

$$\sup_{n \geq 0} \int x^{kn} d\mu = \infty \text{ and } \inf_{n \geq 0} \int x^{kn} d\mu = -\infty.$$

*Proof.* From Corollary 6.4 we know that the convex-hull of the set $\{x^{kn} : n \geq 0\}$ is dense in $C(a, -1)$ with the topology of uniform convergence on compact sets. Thus, since the dual space of $C(a, -1)$ is $M(a, -1)$, the Hahn-Banach Criterion (Theorem 2.2) implies that if $\mu \in M(a, -1)$ is nonzero, then $\sup_{n \geq 0} \int x^{kn} d\mu = \infty$. By applying the previous observation to $(-1)\mu$, we have $\inf_{n \geq 0} \int x^{kn} d\mu = -\infty$. □

## 7. Convex-Cyclic Multiplication Operators

If $\mu$ is a compactly supported positive regular Borel measure on $\mathbb{R}$, then let $L^2_{\mathbb{R}}(\mu)$ denote the real Hilbert space of all real-valued Lebesgue measurable functions that are square integrable with respect to $\mu$. Also let $M_{x,\mu}$ denote the operator of multiplication by the independent variable $x$ acting on $L^2_{\mathbb{R}}(\mu)$. Then $M_{x,\mu}$ is a real symmetric operator on $L^2_{\mathbb{R}}(\mu)$. In this section we determine when $M_{x,\mu}$ is convex-cyclic and characterize its convex-cyclic vectors. We also investigate the invariant convex-sets for $M_{x,\mu}$ in the case where $M_{x,\mu}$ is convex-cyclic.

Recall that a vector $v$ is a *cyclic vector* for a continuous linear operator $T$ on a space $X$ if the *linear span* of the orbit of $v$ under $T$ is dense in $X$. Similarly, a vector $v$ is a *convex-cyclic vector* for $T$ if the *convex-hull* of the orbit of $v$ under $T$ is dense in $X$. Thus, $v$ is a cyclic vector for $T$ if $\{p(T)v : p \text{ is a polynomial}\}$ is dense in $X$ and $v$ is a convex-cyclic vector for $T$ if $\{p(T)v : p \text{ is a convex-polynomial}\}$ is dense in $X$.

It is well known that for a positive compactly supported regular Borel measure $\mu$ on $\mathbb{R}$ that a function $f \in L^2_{\mathbb{R}}(\mu)$ is a cyclic vector for $M_{x,\mu}$ if and only if $|f| > 0$ $\mu$-almost everywhere. This holds simply because the polynomials are dense in the space $C(K)$ of all continuous functions on a compact set $K \subseteq \mathbb{R}$. Hence every invariant subspace for $M_{x,\mu}$ is a zero based invariant subspace. That is, a subspace consisting of all functions that vanish on a given set of positive measure. The following theorem tells us when the multiplication operator $M_{x,\mu}$ is convex-cyclic on $L^2_{\mathbb{R}}(\mu)$ and identifies its convex-cyclic vectors as the same as its cyclic vectors!

**Theorem 7.1.** *If $\mu$ is a positive finite regular Borel measure with compact support in $\mathbb{R}$ and $M_{x,\mu}$ is the operator of multiplication by $x$ on the Hilbert space $L^2_{\mathbb{R}}(\mu)$ of all real-valued functions that are square integrable with respect to $\mu$, then the following statements hold:*

(1) *$M_{x,\mu}$ is convex-cyclic on $L^2_{\mathbb{R}}(\mu)$ if and only if $\mu([-1,\infty)) = 0$;*
(2) *If $M_{x,\mu}$ is convex-cyclic, then the convex-cyclic vectors for $M_{x,\mu}$ are the same as its cyclic vectors;*
(3) *If $M_{x,\mu}$ is convex-cyclic on $L^2_{\mathbb{R}}(\mu)$, then $M^k_{x,\mu}$ is also convex-cyclic for any odd integer $k \geq 1$. Furthermore, $M_{x,\mu}$ and $M^k_{x,\mu}$ have the same convex-cyclic vectors.*

*Proof.* (1) The fact that $\mu([-1,\infty)) = 0$ when $M_{x,\mu}$ is convex-cyclic follow from item (3) of Proposition 2.1. Conversely, if $\mu([-1,\infty)) = 0$, then it follows from Corollary 5.2 that the convex-polynomials are dense in $L^2_{\mathbb{R}}(\mu)$; which implies that $M_{x,\mu}$ is convex-cyclic and the constant function 1 is a convex-cyclic vector for $M_{x,\mu}$.

(2) First recall that the cyclic vectors for $M_{x,\mu}$ are those functions $f \in L^2_{\mathbb{R}}(\mu)$ that satisfy $|f| > 0$ $\mu$ a.e. Now, if $f$ is a convex-cyclic vector for $M_{x,\mu}$, then clearly it is also a cyclic vector and thus must satisfy $|f| > 0$ $\mu$ a.e.. Conversely, suppose that $M_{x,\mu}$ is convex-cyclic and $f \in L^2_{\mathbb{R}}(\mu)$ satisfies $|f| > 0$ $\mu$ a.e. Since $M_{x,\mu}$ is convex-cyclic we know from (1) that $\mu([-1,\infty)) = 0$ and since $|f| > 0$ $\mu$ a.e., Corollary 5.2 item (3) says exactly that $f$ is a convex-cyclic vector for $M_{x,\mu}$.

a

(3) Suppose that $M_{x,\mu}$ is convex-cyclic with convex-cyclic vector $f \in L^2_{\mathbb{R}}(\mu)$ and that $k$ is an odd positive integer. Then by definition $C := \{p(x) \cdot f(x) : p \in \mathcal{CP}\}$ is dense in $L^2_{\mathbb{R}}(\mu)$ and thus $|f| > 0$ $\mu$ a.e. and $\mu([-1, \infty)) = 0$. To show that $\{p(x^k)f(x) : p \in \mathcal{CP}\}$ is dense in $L^2_{\mathbb{R}}(\mu)$ it suffices to show that $\int x^{kn} f(x)g(x)d\mu = \infty$ for every non-zero $g \in L^2_{\mathbb{R}}(\mu)$. So let $g \in L^2_{\mathbb{R}}(\mu)$ be non-zero. Since $|f| > 0$ $\mu$ a.e. and since $g$ is not the zero function, then $fg$ must be non-zero on a set of positive $\mu$ measure. Thus the measure $fg d\mu$ is a non-zero measure carried by the set $(a, -1)$ for some $a < -1$. Since $k$ is an odd positive integer, then it follows from Corollary 6.5 that $\sup_{n \geq 0} \int x^{kn} f(x)g(x)d\mu = \infty$. It now follows from the Hahn-Banach Criterion (Theorem 2.2) that the convex-hull of $\{x^{kn}f(x) : n \geq 0\}$ is dense in $L^2_{\mathbb{R}}(\mu)$. Thus $f$ is a convex-cyclic vector for $M^k_{x,\mu}$. In particular, $M_{x,\mu}$ and $M^k_{x,\mu}$ have the same convex-cyclic vectors (obviously convex-cyclic vectors for $M^k_{x,\mu}$ are also convex-cyclic for $M_{x,\mu}$). □

## 8. Invariant Convex Sets

If $T$ is a continuous linear operator on a locally convex-space $X$ and $C$ is a subset of $X$, then we say that $C$ is invariant under $T$ if $T(C) \subseteq C$. In this section we will characterize the invariant closed convex sets for the operator $M_{x,\mu}$ of multiplication by $x$ on $L^2_{\mathbb{R}}(\mu)$ when $\mu([-1, \infty)) = 0$. In this case there is a surprising answer - the invariant closed convex sets are the same as the invariant closed subspaces!

**Theorem 8.1.** *If $\mu \in M^+_c(\mathbb{R})$, then the following are equivalent:*

(1) *the closed invariant convex sets for the multiplication operator $M_{x,\mu}$ on $L^2_{\mathbb{R}}(\mu)$ are the same as the closed invariant subspaces for $M_{x,\mu}$;*
(2) $\mu([-1, \infty)) = 0$.

*Proof.* (1) $\Rightarrow$ (2) Suppose that every closed invariant convex set for $M_{x,\mu}$ is a subspace. To show $\mu([-1, \infty)) = 0$, suppose that $\mu([-1, \infty)) > 0$. Then either $\mu([-1, 1]) > 0$ or $\mu([0, \infty)) > 0$. In either case we will obtain a contradiction. Let $A = \{f \in L^2_{\mathbb{R}}(\mu) : |f(x)| \leq 1 \text{ for } \mu \text{ a.e. } x \in [-1, 1]\}$ and let $B = \{f \in L^2_{\mathbb{R}}(\mu) : f(x) \geq 0 \text{ for } \mu \text{ a.e. } x \in [0, \infty)\}$. Both $A$ and $B$ are closed convex sets in $L^2(\mu)$ that are invariant under $M_{x,\mu}$. By our assumption that $\mu([-1, \infty)) > 0$, one of $A$ or $B$ is nonzero. Also by assumption both $A$ and $B$ are subspaces. But clearly neither $A$ nor $B$ is a subspace unless they are the zero space. This gives the required contradiction and therefore $\mu([-1, \infty)) = 0$.

(2) $\Rightarrow$ (1) Suppose that $\mu([-1, \infty)) = 0$. We must show that every closed invariant convex set is a subspace. First notice that a convex set is a subspace if and only if it is closed under scalar multiplication. Furthermore for a closed convex set, it suffices to show that it is invariant under multiplication by non-zero scalars. So, let $K$ be a closed convex set in $L^2_{\mathbb{R}}(\mu)$ that is invariant under $M_{x,\mu}$, let $f \in K \setminus \{0\}$, and let $c \in \mathbb{R} \setminus \{0\}$. Consider the measure $\nu = |f|^2 d\mu$ and the

space $L^2_\mathbb{R}(|f|^2 d\mu)$. Since $\mu([-1,\infty)) = 0$, it follows that $\nu([-1,\infty)) = 0$ and, by Corollary 5.2, we know that the convex-polynomials are dense in $L^2(\nu)$. Thus, there are non-zero convex-polynomials $\{p_n\}_{n=1}^\infty$ such that $p_n \to c$ in $L^2(\nu)$. Hence $\int |p_n f - cf|^2 d\mu = \int |p_n - c|^2 \cdot |f|^2 d\mu = \int |p_n - c|^2 d\nu \to 0$ as $n \to \infty$ and so $p_n f \to cf$ in $L^2(\mu)$. Now since $K$ is convex and invariant under $M_{x,\mu}$ and $f \in K$, $p_n f \in K$ for every $n$ and $p_n f \to cf$ in $L^2(\mu)$. Since $K$ is a closed set, $cf \in K$. It now follows that $K$ is invariant under scalar multiplication. So $K$ is a subspace, as desired. □

## 9. Approximation on the Positive Real Line

We have shown that convex-polynomials can be used to approximate various classes of functions defined on bounded subsets of $(-\infty, -1)$. One may ask what functions can be approximated by convex-polynomials on subsets of $[-1, \infty)$. In this section we show that such functions must have a convex-power series representation.

A power series is a *convex-power series* if it is an infinite convex combination of the monomials $\{1, x, x^2, \ldots\}$. That is, if it has the form $\sum_{n=0}^\infty c_n x^n$ where $c_n \geq 0$ for all $n \geq 0$ and $\sum_{n=0}^\infty c_n = 1$. Note that since the sum of the coefficients converges that a convex-power series will converge absolutely on the interval $[-1, 1]$. In fact, the series converges absolutely and uniformly on the closed disk $\{z \in \mathbb{C} : |z| \leq 1\}$.

**Example 9.1.** Examples of functions with convex-power series representations include
$$e^{(x-1)} = \sum_{n=0}^\infty \frac{e^{-1}}{n!} x^n \quad \text{and} \quad \frac{1-a}{(x-a)} = \sum_{n=0}^\infty \frac{(a-1)}{a^{n+1}} x^n,$$
where $a > 1$ as well as products, compositions, and convex-combinations of these functions.

We will say that a real-valued function $f$ defined on an interval $(a, b)$ has a *convex-power series representation* if $f$ is the restriction of a convergent convex-power series to the interval $(a, b)$. That is, if there exists a convex-power series $\sum_{n=0}^\infty c_n x^n$ that converges on an interval $(-R, R)$ which contains the interval $(a, b)$, and satisfies that for every $x \in (a, b)$ we have $f(x) = \sum_{n=0}^\infty c_n x^n$.

**Proposition 9.2.** *(a) If $\mathcal{C}$ is a set of convex-polynomials and there exists a $c > 0$ such that the set $\{p(c) : p \in \mathcal{C}\}$ is bounded, then the set $\mathcal{C}$ is a normal family on the disk $\{z \in \mathbb{C} : |z| < c\}$.*
*(b) The set $\mathcal{CP}$ of all convex-polynomials is a normal family on $\mathbb{D}$.*

*Proof.* (a) If $p$ is a convex-polynomial, then its coefficients are all non-negative and the triangle inequality implies $|p(z)| \leq p(|z|)$. Since a convex-polynomial is increasing on the interval $[0, \infty)$, we obtain the following: if $p \in \mathcal{C}$, then $|p(z)| \leq p(|z|) \leq p(c)$ for all $z$ satisfying $|z| \leq c$. Since $\mathcal{C}$ is a uniformly bounded family of analytic functions, the result follows from Montel's theorem.

(b) This follows from (a) with $c = 1$ and noting that $p(1) = 1$ for every convex polynomial. □

**Theorem 9.3.** *If $f$ is a real valued function defined on a closed and bounded interval $I = [a, b]$ that is contained in $[-1, \infty)$, then the following are equivalent:*

  (1) *$f$ is the pointwise limit of a sequence of convex-polynomials on $I$.*
  (2) *$f$ is the uniform limit of a sequence of convex-polynomials on $I$.*
  (3) *$f$ has a convex-power series representation on $I$ which converges absolutely and uniformly on compact subsets of $(-R, R)$ where $R \geq \max\{1, |a|, |b|\}$.*

*Proof.* Clearly, (3) $\Rightarrow$ (2) $\Rightarrow$ (1). We will prove that (1) $\Rightarrow$ (3). So, assume that there exists a sequence $\{p_n\}_{n=1}^{\infty}$ of convex-polynomials such that $p_n(x) \to f(x)$ for every $x \in I$ where $I = [a, b] \subseteq [-1, \infty)$. Let $R = \max\{1, |a|, |b|\}$. By Proposition 9.2 the set $\{p_n\}$ is a normal family on $B = \{z : |z| < R\}$ and so a subsequence $\{p_{n_k}\}$ of $\{p_n\}$ converges uniformly on compact subsets of $B$ to an analytic function $g$ which has a convex-power series representation on $B$. Since $\{p_{n_k}\}$ also converges to $f$ on $[a, b]$, $f$ has a convex-power series representation on $[a, b]$. □

## 10. A Question

The focus of this paper was on approximation by convex-polynomials. A convex-polynomial $p$ is distinguished by it's behavior at the two points $x = 0$ and $x = 1$. Namely, a polynomial $p$ is convex if and only if $p^{(k)}(0) \geq 0$ for all $k \geq 0$ and $p(1) = 1$. It is natural to ask: which functions can be approximated by polynomials whose values and derivatives are restricted at a finite number of points?

Dept. of Mathematics, Washington and Lee University, Lexington VA 24450
*E-mail address*: `feldmanN@wlu.edu`

Dept. of Mathematics, Bucknell University, Lewisburg PA 17837
*E-mail address*: `pmcguire@bucknell.edu`